\documentstyle[12pt,amssymb]{article}
\input xy
\xyoption{all} \hfuzz=1.5pt \tolerance=500

\newcommand{\ot}{\otimes}

\newcommand{\tr}{\triangleright}
\newcommand{\tl}{\triangleleft}

\newcommand{\va}{\varphi}
\newcommand{\Om}{\Omega}

\newcommand{\id}{{\bf 1}}
\newcommand{\co}{{\B C}}
\newcommand{\s}{\sigma}
\newcommand{\lm}{\lambda}

\newcommand{\bd}{\begin{document}}
\newcommand{\ed}{\end{document}}
\newcommand{\Htp}{\mathbin{\stackrel{\cdot}{\bigotimes}}}
\newcommand{\Hts}{\mathbin{\stackrel{\cdot}{\bigoplus}}}
\newcommand{\nl}{\nu\in\Lambda}
\newcommand{\xl}{X_\nu;\nl}
\newcommand{\x}{X_\nu}
\newcommand{\la}{\langle}
\newcommand{\ra}{\rangle}
\newcommand{\ch}{\cal H}
\newcommand{\su}{\subseteq}
\newcommand{\B}{\Bbb}
\newcommand{\e}{\varepsilon}
\newcommand{\vk}{\varkappa}
\newcommand{\vt}{\vartheta}
\newcommand{\al}{\alpha}
\newcommand{\de}{\delta}
\newcommand{\lo}{\Longleftrightarrow}
\newcommand{\k}{{\cal F}}
\newcommand{\cK}{{\cal F}}
\newcommand{\bb}{{\cal B}}
\newcommand{\dd}{{\cal D}}
\newcommand{\r}{{\cal R}}
\newcommand{\cR}{{\cal R}}
\newcommand{\f}{{\cal K}}
\newcommand{\cb}{{\cal CB}}
\newcommand{\wrr}{\widetilde{\cal R}}
\newcommand{\br}{^\bullet{\cal R}}
\newcommand{\en}{E_\nu}
\newcommand{\el}{{\cal L}}
\newcommand{\lr}{\Longrightarrow}
\newcommand{\Long}{\Longleftarrow}
\newcommand{\q}{\quad}
\newcommand{\qq}{\qquad}
\newcommand{\cd}{\cdot}
\newcommand{\fur}{f: E\to{\B R}}
\newcommand{\furo}{f_0: E_0\to{\B R}}
\newcommand{\fuc}{f: E\to{\B C}}
\newcommand{\ii}{\infty}
\newcommand{\di}{\diamondsuit}
\newcommand{\bgd}{{\bigtriangledown}}
\newcommand{\bu}{{\bigtriangleup}}
\newcommand{\bc}{{completely bounded}}
\newcommand{\cc}{{completely contractive}}
\newcommand{\qs}{{quantum space}}
\newcommand{\isc}{{isometric}}
\newcommand{\ism}{{isomorphism}}
\newcommand{\qss}{{quantum spaces}}
\newcommand{\bco}{{completely bounded operator}}
\newcommand{\bcos}{{completely bounded operators}}
\newcommand{\res}{{respectively}}
\newcommand{\tp}{{tensor product}}
\newcommand{\eq}{{equivalent}}
\newcommand{\qtp}{{quantum tensor product}}
\newcommand{\mma}{\mathrel{\mathop{\otimes}\limits_{A}}}
\newcommand{\mmg}{\mathrel{\mathop{\otimes}\limits_{{c_0}}}}
\newcommand{\mmb}{\mathrel{\mathop{\otimes}\limits_{\bb}}}
\newcommand{\mmp}{\mathrel{\mathop{\otimes}\limits_{p}}}
\newcommand{\mmh}{\mathrel{\mathop{\otimes}\limits_{h}}}
\newcommand{\mmf}{\mathrel{\mathop{\otimes}\limits_{4}}}
\newcommand{\mmi}{\mathrel{\mathop{\otimes}\limits_{i}}}
\newcommand{\mmm}{\mathrel{\mathop{\otimes}\limits_{\bb-\bb}}}
\newcommand{\mmo}{\mathrel{\mathop{\cdot}\limits_{1}}}
\newcommand{\mmA}{\mathrel{\mathop{\otimes}\limits_{A-A}}}
\newcommand{\mmd}{\mathrel{\mathop{\cdot}\limits_{2}}}
\newcommand{\msp}{\mathrel{\mathop{\otimes}\limits_{sp}}}
\newcommand{\mms}{\stackrel{h}{\otimes}}
\newcommand{\mmt}{\stackrel{4}{\otimes}}
\newcommand{\mmx}{\stackrel{i}{\otimes}}
\newcommand{\mmy}{\stackrel{p}{\otimes}}
\newcommand{\mmz}{\stackrel{sp}{\otimes}}
\newcommand{\gd}{\ddagger}
\newcommand{\od}{\odot}
\newcommand{\mt}{\mapsto}
\newcommand{\mmc}{\mathrel{\mathop{\otimes}\limits_{\sim}}}
\newcommand{\mme}{\stackrel{\sim}{\otimes}}

\begin{document}

   \centerline{{ \bf Extreme flatness and Hahn-Banach type theorems  }}

   \centerline{{\bf for normed modules over $c_0$}\footnote{This research
was supported by the Russian Foundation for Basic Research (grant No. 08-01-00867).}} 

   \vspace{1cm}

\centerline{A.~Ya.~Helemskii} 

\centerline{Faculty of Mechanics and Mathematics} 

\centerline{Moscow State University} 

\centerline{Moscow 119992 Russia} 

\bigskip

\bigskip

\centerline{\bf Introduction: formulation of the main results and comments } 

\bigskip

\noindent In this paper we consider a certain specific case of a well known typical question in the 
theory of normed algebras and their modules. This is a question about conditions that ensure the 
preservation of isometries under projective tensor product of modules. Such a question is 
intimately connected with the problem of extension of a given bounded morphism from a submodule to 
a bigger module with the preservation of its norm. In other words, it is connected with the 
question of the existence, in certain situations, of module versions of the classical Hahn-Banach 
Theorem. 

We proceed to relevant formal definitions. 

Let $A$ be a normed algebra. We shall use the symbol ` $\mma$ ' for the {\it non-completed} 
projective module tensor product of $A$-modules and of their
bounded morphisms. (See, e.g.,\cite{he5} or, as to the initial `completed' version, the pioneering 
paper of Rieffel~\cite{rie1} or the textbooks~\cite[II.3]{he2}~\cite[VI.3.2]{he3}). 

The identity operator on a linear space (or a module) $Z$ will be denoted by $\id_Z$, or just 
$\id$, if there is no danger of misunderstanding. 

Let us distinguish a class, so far arbitrary, of right normed $A$-modules and denote it by ${\cal 
K}$. In the spirit of the well-known definitions of a flat and of a strictly flat Banach module 
(\cite[VII.1]{he2},~\cite [VII.1.3]{he3}), we give the following

\medskip
{\bf Definition.} A normed left $A$-module $Z$ is called {\it extremely flat with respect to the 
class ${\cal K}$} or, for short, {\it ${\cal K}$-E-flat}, if, for every isometric morphism 
 ${\textsf i}:X\to Y$ of right modules, belonging to ${\cal K}$, the operator ${\textsf i}\mma\id_Z:
 X\mma Z\to Y\mma Z$ is also isometric.
 
 \medskip 
If $A:=\co$, that is if we deal with just normed spaces, the well known theorem of 
Grothendieck~\cite[Thm. 1]{gro}, being adapted to non-complete spaces, gives a full description of 
the extremely flat objects in the following way. 
{\it A normed space} (`normed $\co$-module') {\it is extremely flat with respect to the class of 
all normed spaces if and only if it is isometrically isomorphic to a dense subspace of 
$L_1(\Omega,\mu)$ for some measure space $(\Omega,\mu)$.} 


\medskip
{\bf Definition.} A normed right $A$-module $Z$ is called {\it extremely injective with respect to 
the class ${\cal K}$} or, for short, {\it ${\cal K}$-E-injective}, if, for every isometric morphism 
 ${\textsf i}:X\to Y$ of right modules, belonging to ${\cal K}$, and for every bounded morphism 
 $\va:X\to Z$ of right $A$-modules, there exists a bounded 
morphism of right modules $\psi:Y\to Z$ such that the diagram 
$$
\xymatrix{ X \ar[r]^{{\textsf i}} \ar [d]^\va
& Y \ar[dl]^\psi\\
 Z }
$$
\noindent  is commutative  and $\|\va\|=\|\psi\|$. In other words, every bounded morphism of right 
modules from $X$ into $Z$ can be extended, after the identification of $X$ with a submodule of $Y$, 
to a morphism from $Y$ to $Z$ with the same norm. 

Thus the assertion that a certain module $Z$ is ${\cal K}$-E-injective can be considered  as a 
`Hahn-Banach type' theorem for given $A$ and ${\cal K}$, with $Z$ playing the role of $\co$ in the 
mentioned classiacal theorem. 

\medskip
If again $A:=\co$, then the extremely injective objects are described by a theorem, connected with 
the names of Nachbin, Goodner, Hasumi and Kelley (see~\cite[p.123]{day} or~\cite[Thm. 
25.5.1]{sem}), which can be easily adapted  to the non-complete case. Namely, {\it a normed space} 
{\it is extremely injective with respect to the class of all normed spaces if and only if it is  
isometrically isomorphic to the space $C(\Om)$, where $\Om$ is an extremely disconnected compact 
space.} 

\medskip
{\bf Remark.} The word  `extremely' in both definitions is chosen because isometric operators or 
morphisms are exactly the so-called extreme monomorphisms in some principal categories of spaces or 
modules in functional analysis (cf., e.g.,~\cite[p. 4]{clm},~\cite[Ch. 0.5]{he1}). 

\medskip
The both introduced notions are closely connected. The link is provided by a proper 
functional-analytic version of the algebraic `law of adjoint associativity'. This version was 
established by Rieffel~\cite{rie2} (who considered Banach modules). With its help, we shall prove 
below (see the beginning of Section 2) the following easy 

\medskip
{\bf Proposition}. {\it Let $A, {\cal K}$ and $Z$ be as above. Then $Z$ is ${\cal K}$-E-flat if and 
only if its dual normed left $A$-module $Z^*$ is ${\cal K}$-E-injective}. 

\bigskip
The notions, defined above, were actually introduced in~\cite{he4}, however, for only some special 
algebras and modules. Namely, the role of a base algebra was played by $\bb(H)$ for a Hilbert space 
$H$, and the class ${\cal K}$ consisted of the so-called semi-Ruan $\bb(H)$-modules. (Speaking 
informally, these are modules, satisfying a proper one-sided version of Ruan axioms for an operator 
space; cf.~\cite{he5}). It was shown that certain $\bb(H)$-modules are extremely flat with respect 
to that ${\cal K}$, and certain Hahn-Banach type theorems for modules over $\bb(H)$ were obtained 
as corollaries. These theorems, in their turn, led to a transparent new proof of one of basic 
theorems of operator space theory, the Arveson-Wittstock Theorem about extensions of completely 
bounded operators (see, e.g.,\cite{efr} or~\cite{he5}). 

Afterwards the results of~\cite{he4} were generalized and considerably strengthened by 
Wittstock~\cite{wit}, who, in particular, replaced $\bb(H)$ by an arbitrary properly infinite 
$C^*$-algebra and established that every semi-Ruan module is ${\cal K}$-E-flat. As an application 
of his results, Wittstock presented a new transparent proof of the 
 Arveson-Wittstock Theorem in a more sophisticated version, that for operator modules.

\bigskip
After the cited papers it seemed natural to look for extremely flat modules over other classes of 
normed algebras and, accordingly, for related Hahn-Banach type theorems. In particular, what can we 
find, if we turn to commutative algebras ? This class, in a sense, is opposite to highly 
non-commutative algebras of~\cite{he4} and~\cite{wit}.

In the present paper we exclusively deal with the apparently simplest of all infinite-dimensional 
commutative normed algebras. This is the algebra $c_0$ of complex-valued sequences, converging to 
0, with the coordinate-wise operations and the uniform norm. It turned out that even in this case 
there is something to say. (Speaking very roughly, extremely flat $c_0$-modules form much larger 
family that one could initially expect). 

\bigskip 
We recall that a normed module $X$ over a normed algebra $A$ is called contractive, if we have 
$\|a\cd x\|\le \|a\|\|x\|$, or, accordingly, $\|x\cd a\|\le \|a\|\|x\|$ for all $a\in A, x\in X$. 
{\it Throughout this paper, all normed modules are always supposed to be contractive}. 

\smallskip
If $A$ and $X$ are as before, we denote the closure of the linear span of the set $\{a\cd x: a\in 
A, x\in X\}$ by $X_{es}$ and call it {\it essential part of $X$}. It is, of course, a submodule. A 
left $A$-module $X$ is called {\it essential} (they often say also `non-degenerate'), if we have 
$X=X_{es}$. The quotient normed $A$-module $X/X_{es}$ is denoted by $X_{an}$; obviously it has zero 
outer multiplication. The {\it annihilator of $A$ in $X$} is the closed left submodule $\{x: a\cd 
x=0$ for all $a\in A\}$ in $X$, denoted by $Ann X$. The quotient left normed $A$-module $X/Ann X$ 
is called the {\it reduced module of $X$} and denoted by $X^{red}$. 

As usual, we call a left $A$-module $X$ {\it faithful}, if $Ann X=0$. Of course, the reduced module 
of every module is faithful. It is easy to show that every essential left $A$-module is faithful 
provided $A$ has a bounded left approximate identity. 

Recall what happens if $A$ is commutative, as it is the case with $c_0$. Then every left $A$-module 
is a right $A$-module with the same bilinear operator of the outer multiplication, and vice versa. 
Therefore we identify both types of modules and say just `$A$-module'. Accordingly, we can speak 
about module projective tensor product of two normed $A$-modules and of two bounded morphisms of 
normed $A$-modules. Moreover, we immediately see that the mentioned tensor product of two modules, 
say $X$ and $Y$, is itself a normed contractive $A$-module with the outer multiplication, well 
defined by $a\cd(x\mma y):=(a\cd x)\mma y$ (or $:=x\mma(a\cd y)$). Besides, the mentioned tensor 
product of two bounded morphisms of normed $A$-modules is obviously itself a bounded morphism of 
the respective modules. 

\bigskip
The main result of the paper gives, within a certain reasonable class of normed $c_0$-modules, a 
full description of extremely flat modules with respect to that class. After some preliminary note, 
we proceed to the definition this class. 

One can immediately see, what makes the work with $c_0$ easier than with other algebras. It is the 
presence in this space of a distinguished countable Schauder basis, consisting of irreducible 
idempotent generators. We mean, of course, the `orts'  $(0,\dots,0,1,0,0\dots)\in c_0$. The $n$-th 
ort (that with 1 as its $n$-th term) will be denoted by ${\bf p}^n$. If $X$ is a normed 
$c_0$-module, we set $X_n:=\{{\bf p}^n\cd x; x\in X\}$ for every $n=1,2,\dots$. We see that $X_n$ 
is a submodule of $X$; it will be called the {\it $n$-th coordinate submodule}. Often, when there 
is no danger of confusion, for $x\in X$ we shall write $x_n$ instead of ${\bf p}^n\cd x$. Of 
course, we have ${\bf p}^n\cd x_n=x_n$. 

\medskip 
{\bf Definition}. A $c_0$-module $X$ is called {\it homogeneous} if, for every $x,y\in X$, the 
equalities $\|x_n\|=\|y_n\|$, for all $n$, imply that $\|x\|=\|y\|$. 

\medskip 
In particular, all essential normed $c_0$-modules, consisting of complex-valued sequences, are  
homogeneous (Proposition 3.1 below). Besides, $l_p$-sums; $1\le p\le\ii$ of arbitrary families of 
normed spaces are obviously homogeneous. (In both cases we mean the coordinate-wise outer 
multiplication). 

It is evident that every homogeneous normed $c_0$-module is faithful.

In this paper, by ${\cal H}$ we denote the class of all homogeneous normed $c_0$-modules, and by 
${\cal H}_{es}$ its subclass, consisting of essential modules. 
  

\medskip 
{\bf Theorem I}. {\it Let $Z$ be an essential (respectively, arbitrary) homogeneous normed 
$c_0$-module. Then $Z$ is extremely flat with respect to ${\cal H}$ (respectively, with respect to 
${\cal H}_{es}$) if and only if, 
for every n, its $n$-th coordinate submodule is isometrically isomorphic to a dense subspace of the 
space $L_1(\Om_n,\mu_n)$ for some measure space $(\Om_n,\mu_n)$.} 

\medskip 
Note that `only if' part of this theorem relies heavily on the theorem of Grothendieck, cited 
above, and it is rather easy corollary of the latter. As to the `if' part, our proof of this is 
more complicated, and it does not use the Grothendieck Theorem). 

\medskip
In fact, we shall prove this theorem in a slightly stronger form; see Proposition 3.3 and Theorem 
3.7 below. 

The following theorem is a rather easy corollary of Theorem I.

\medskip 
{\bf Theorem II} (see end of Section 4). {\it Let $Z$ be an essential (respectively, arbitrary) 
homogeneous normed $c_0$-module. Then the dual module $Z^*$ is extremely injective with respect to 
${\cal H}$ (respectively, with respect to ${\cal H}_{es}$)  if and only if for every n we have that 
its $n$-th coordinate submodule $(Z^*)_n$ is isometrically isomorphic to the Banach space 
$L_\ii(\Om_n,\mu_n)$ for some measure space $(\Om_n,\mu_n)$.} 

\medskip 
In particular, all $c_0$-modules $l_p; 1\le p<\ii$ are ${\cal H}$-E-flat whereas the same $l_p$ and 
also $l_\ii$ are ${\cal H}$-E-injective.

\medskip 
In both theorems we assumed that some participating modules are essential. Such a condition can not 
be omitted: a non-essential homogeneous normed module (being always ${\cal H}_{es}$-E-flat) is not 
bound to be ${\cal H}$-E-flat. As a matter of fact, {\it the $c_0$-module $l_\ii$} (apparently the 
first faithful non-essential $c_0$-module that comes in mind), {\it is not extremely flat with 
respect to the class of all homogeneous modules.} This is Theorem 4.3. 

\bigskip
Let us make some comments on the proof of the main result. In the very beginning we observe that, 
under some conditions, tensor products of $c_0$-modules and their morphisms can be described in a 
rather transparent and `workable' form (Proposition 1.6). In particular, this is helpful in making 
the principal preparatory step, Lemma 3.4 of somewhat technical character. At the end of our 
argument, we have used the following fact: if $X$ or $Z$ are essential, then 
the property of $\va:X\to Y$ to be (just) injective implies the same property of $\va\mmg\id_Z$. 

\bigskip
Thus, trying to prove the preservation of isometries, we came across another typical question of 
the theory of normed algebras. Which conditions ensure the preservation, under projective tensor 
multiplication of modules, of the property of a given morphism to be injective ? We believe that 
such a question deserves to be considered independently. Of course, it sounds similar to its well 
known pure algebraic prototype, which leads to the fundamental notion of the (algebraic) flatness. 
But here we deal with the bounded morphisms and a kind of functional-analytic tensor product. This 
profoundly affects the situation.  

As a matter of fact (see Example 2.3), if $X,Y,Z$ are normed $c_0$-modules, even consisting of 
sequences, then it can well be that a bounded morphism $\va:X\to Y$ is injective whereas 
$\va\mmg\id:X\mmg Z\to Y\mmg Z$ is not. However, if we are given arbitrary  normed $c_0$-modules  
$X,Y,Z$ and a {\it topologically injective} (in particular, isometric) 
 morphism $\va:X\to Y$ then  {\it the operator 
$\va\mmg\id:X\mmg Z\to Y\mmg Z$ is also injective}. (Note that at the same time
it is not bound to be again topologically injective). This is the future Theorem 2.4. 

\medskip 
{\bf Remark.} We want to emphasize that we work in this paper, in a similar way as 
in~\cite{he4}\cite{wit}, with the non-completed version of the module projective tensor product. If 
we replace the latter by the respective completed version, Theorem 2.4 fails to be true. One can 
easily construct respective counter-examples, taking some spaces without the approximation 
property. 

\bigskip
{\centerline {\bf 1. Some preparations }} 

\bigskip
We begin our preliminaries with a proposition of somewhat general character. In particular, it will 
enable us to derive Theorem II from Theorem I. This proposition actually appeared in~\cite[Prop. 
9]{he4}, but in a certain special case and in a slightly disguised form. 

In what follows $A$ is a normed algebra, so far arbitrary, and ${\bf h}_A(\cd,\cd)$ is the symbol 
of the space of all bounded morphisms between right normed  modules. Such spaces are equipped with 
the operator norm. 

\medskip
{\bf Proposition 1.1}. {\it Let $X$ and $Y$ be right normed $A$-modules, $Z$ a left normed 
$A$-module, ${\textsf i}:X\to Y$ an isometric morphism, $Z^*$  the right Banach $A$-module, dual to 
$Z$. Then the following statements are equivalent: 

(i) the operator ${\textsf i}\mma\id_Z:X\mma Z\to Y\mma Z$ is an isometry

(ii) for every bounded morphism $\va:X\to Z^*$ of right $A$-modules, there exists a bounded 
morphism of right modules $\psi:Y\to Z^*$ such that the diagram 
$$
\xymatrix{ X \ar[r]^{{\textsf i}} \ar [d]^\va
& Y \ar[dl]^\psi\\
 Z^* }
$$
\noindent  is commutative  and $\|\va\|=\|\psi\|$. }

\smallskip
$\tl$ According to the functional-analytic version of the law of the adjoint associativity 
(cf.~\cite{rie2} or~\cite[Ch. 8.0]{he5}) the normed space {\bf h}$_A(X,Z^*)$ coincides with the 
space $(X\mma Z)^*$ up to the isometric isomorphism, taking a morphism $\va:X\to Z^*$ to the 
functional $f:X\mmb Z\to\co$, well-defined by $f(x\mma z)=[\va(x)](z)$. Similarly, {\bf 
h}$_A(Y,Z^*)$ is identified with $(Y\mma Z)^*$. Moreover, one can easily check that we have a 
commutative diagram 
$$
\xymatrix@C+20pt{ {\bf h}_A(Y,Z^*) \ar[r]^{{\textsf i}_*} \ar[d]
& {\bf h}_A(X,Z^*) \ar[d]\\
(X\mma Z)^* \ar[r]^{{\textsf i}^\bullet} & {(Y\mma Z)^*} }.
$$
\noindent Here the vertical arrows depict isometric isomorphisms of normed spaces, acting as it was 
indicated, ${\textsf i}_*$ acts as $\beta\mt\beta{\textsf i}$, and ${\textsf i}^\bullet$ is the 
operator which is adjoint to ${\textsf i}\mma\id_Z:X\mma Z\to Y\mma Z$. 

It is obvious that the assertion (ii) is equivalent to the following statement: the operator 
${\textsf i}_*$ maps the closed unit ball in the domain space onto the closed unit ball in the 
range space. Because of the diagram above, this assertion, in its turn, is equivalent to the 
statement that ${\textsf i}^\bullet$ has the same property. But, as an obvious corollary (in fact, 
an equivalent formulation) of the Hahn-Banach theorem, an adjoint operator has the indicated 
property if and only if the original operator is {\isc}. The rest is clear. $\tr$ 

\medskip
An immediate corollary is Proposition that was formulated in the beginning of Introduction.

\medskip
As a byproduct, we have the following observation. 

\medskip
{\bf Proposition 1.2.} {\it  Suppose that $X,Y,Z$ and ${\textsf i}$ are as before, and $Z_0$ is a 
dense submodule of $Z$. Then ${\textsf i}\mma\id_Z$ is an isometry  if and only if the same is true 
of  ${\textsf i}\mma\id_{Z_0}$.} 

\smallskip
$\tl$ Indeed, the dual modules of $Z$ and $Z_0$ coincide, and therefore the assertion (ii) above is 
valid  if and only if it is valid after the replacing of $Z$ by $Z_0$. The rest is clear. $\tr$ 

\medskip
Later we shall come across quite a few diagrams like that one above. To write down them all would 
take too much space. In this connection the following terminology is convenient. We shall say that 
the morphisms $\va:X_1\to X_2$ and $\psi:Y_1\to Y_2$ acting between normed $A$ -modules, are {\it 
isometrically equivalent}, if there exist {\isc} {\ism}s of $A$-modules $I$ and $J$ such that the 
diagram 
$$
\xymatrix@C+20pt{ X_1 \ar[r]^{\va} \ar[d]_{I}
& X_2 \ar[d]^{J}\\
Y_1 \ar[r]^{\psi} & {Y_2} }\eqno(1.1)
$$
\noindent is commutative. In particular, we shall speak about the isometric equivalence of two 
operators ($\co$-modules). As to the isomorphisms $I$ and $J$, we shall say that they {\it 
implement} the mentioned kind of the equivalence. 

\bigskip
From now on we concentrate on the case $A:=c_0$. We  need some further notation and several 
elementary facts, concerning $c_0$-modules and their tensor products. 

Let $X$ be an arbitrary $c_0$-module, $X_n; n=1,2,...$ its coordinate submodules (see 
Introduction). Note that the outer multiplication in $X_n$ acts as $\xi\cd x=\xi_nx$. We denote by 
$\al_n^X:X_n\to X$ the  respective natural embeddings, and by \\ $\beta_n^X:X\to X_n$ the 
projections $x\mt x_n$. Clearly, we have morphisms of $c_0$-modules that are isometries and, 
respectively, coisometries ( = quotient maps). 

 For every $N=1,2,\dots$ we set $P^N:=\sum_{n=1}^N{\bf p}^n\in c_0$. It is easy to see that {\it for every 
 $x\in X_{es}$} (cf. Introduction)  {\it  we have} 
$$
x=\lim_{N\to\ii}P^N\cd x.\eqno(1.2)
$$

Consider the pure algebraic $c_0$-module ${\textsf X }_{n=1}^\ii X_n$, consisting of all sequences 
$(x_1,...,x_n,...);x_n\in X_n$ and endowed with the coordinate-wise operations. Introduce the map 
$$
\s^X:X\to{\textsf X}_{n=1}^\ii X_n:x\mt(x_1,\dots,x_n,\dots);
$$
 this 
is, of course, a $c_0$-module morphism. Obviously, $Ker(\s^X)$ coincides with $Ann(X)$, and hence 
it is closed. Therefore we can (and will) identify the submodule $Im(\s^X)$ in ${\textsf X 
}_{n=1}^\ii X_n$ with $X^{red}$ (cf. Introduction) and endow it with the respective quotient norm. 

We see that $\s^X$ is injective if and only if $X$ is faithful. In particular, if $X$ is essential
and a fortiori faithful (see (1.2)), $\s^X$ is certainly injective. 

If $x\in X_n$, then the sequence $\s^X(x)=(0,\dots,0,x,0,\dots)$ belongs to $(X^{red})_n$. Taking 
into account that $\|y\|\ge\|x\|$ for all $y$ with $\s^X(y)=\s^X(x)$, we immediately obtain 

\medskip
{\bf Proposition 1.3}. {\it The birestriction $\s^X_n:X_n\to (X^{red})_n$ of  $\s^X$ is an 
isometric isomorphism.} $\tl\;\tr$ 

\medskip
{\bf Proposition 1.4}. {\it For every $n$, the $c_0$-modules $(X^*)_n$ and $(X_n)^*$ are 
isometrically isomorphic.} 

\smallskip
$\tl$ Morphisms $(\al_n^X)^*\al_n^{X^*}:(X^*)_n\to (X_n)^*$ and 
$\beta_n^{X^*}(\beta_n^X)^*:(X_n)^*\to (X^*)_n$ are contractive and inverse to each other. $\tr$

\bigskip 
Now let $Z$ be another $c_0$-module. Our object of interest is the $c_0$-module $X\mmg Z$ (cf. 
Introduction). 

Throughout the paper, $\mmp$ will be the symbol of the non-completed projective tensor product of 
normed spaces (= $\co$-modules). The projective tensor norm will be denoted by $\|\cd\|_p$.

\medskip
{\bf Proposition 1.5}. {\it There exists an isometric isomorphism $\rho_n^{X,Z}:X_n\mmp Z_n\to 
(X\mmg Z)_n$, well defined by $x\ot z\mt x\mmg z$.} 

\smallskip
$\tl$ Consider the contractive linear operators $\rho:X_n\mmp Z_n\to X\mmg Z$ and \\ $\pi:X\mmg 
Z\to X_n\mmp Z_n$, associated with the contractive bilinear operator 

\smallskip
$X_n\times Z_n\to X\mmg Z:(x,z)\mt x\mmg z$ and the contractive balanced bilinear operator $X\times 
Z\to X_n\mmp Z_n:(x,z)\mt {\bf p}^n\cd x\ot {\bf p}^n\cd z$, respectively. Since $\pi\rho=\id$, we 
obtain that $\rho$ is an isometry (whereas $\pi$ is a coisometry). Obviously, 
the image of $\rho$ is exactly $(X\mmg Z)_n$. It remains to denote by $\rho_n^{X,Z}$ the respective 
corestriction. $\tr$ 

\medskip
Now we turn to the normed module $(X\mmg Z)^{red}$ and to the coisometric morphism $\s^{X,Z}:X\mmg 
Z\to(X\mmg Z)^{red}$, which is, by definition, the respective corestriction of $\s^{X\mmg Z}$ (cf. 
above).
We want to describe them, up to an isometric isomorphism and, respectively, isometric equivalence, 
in terms, convenient for their study. 

\smallskip
Consider the pure algebraic $c_0$-module ${\textsf X }_{n=1}^\ii(X_n\ot Z_n)$ with the 
coordinate-wise operations. For $x\in X,z\in Z$ we shall denote by $x\odot z$ the sequence \\
$x_1\ot z_1,\dots,x_n\ot z_n,\dots)$, belonging to this module. Denote by $X\odot Z$ the submodule 
of ${\textsf X }_{n=1}^\ii(X_n\ot Z_n)$, defined as the linear span of all such sequences. 

Introduce a bilinear operator $X\times Z\to X\odot Z:(x,z)\mt x\odot z$; clearly it is balanced. 
Therefore it gives rise to the linear operator and, obviously, a surjective $c_0$-module morphism 
$\odot_{X,Z}:X\mmg Z\to X\odot Z$, well defined by $x\mmg z\mt x\odot z$. 

For $v\in X\odot Z$ we set 
$$
\|v\|_\odot:=\inf\{\sum_{k=1}^m\|x^k\|\|z^k\|\},\eqno(1.3)
$$
 where the infimum is taken over all representations of $v$ in the form $\sum_{k=1}^mx^k\odot z^k; x^k\in X,
 z^k\in Z$. 

\medskip
{\bf Proposition 1.6.} {\it The function $v\mt\|v\|_\odot$ is a norm on $X\odot Z$. Moreover, with 
respect to this norm $X\odot Z$ is isometrically isomorphic to $(X\mmg Z)^{red}$, and $\odot_{X,Z}$ 
is isometrically equivalent to $\s^{X,Z}$. In more details, there is a commutative diagram 
$$
\xymatrix@C+20pt{X\mmg Z \ar[r]^{\s^{X,Z}} \ar[d]_{\id}
&(X\mmg Z)^{red} \ar[d]^{\iota^{X,Z}}\\
X\mmg Z \ar[r]^{\odot_{X,Z}} & {X\odot Z} }\eqno(1.4)
$$
where $\iota^{X,Z}$ is an isometric isomorphism of $c_0$-modules.} 

\smallskip
$\tl$ Since  $\odot_{X,Z}$ is surjective, $X\odot Z$ is a seminormed module with respect to the 
seminorm $\|v\|':=\inf\{\|u\|;\odot_{X,Z}(u)=v\}$. 

\smallskip
First, we shall show that $\|\cd\|_\odot=\|\cd\|'$. Indeed, taking an arbitrary representation 
$v=\sum_{k=1}^mx^k\odot z^k$ and looking at $u:=\sum_{k=1}^mx^k\mmg z^k\in X\mmg Z$, we easily see 
that  $\|v\|'\le\|v\|_\odot$. 
On the other hand, for every $\e>0$ we can take $u\in X\mmg Z$ with $\odot_{X,Z}(u)=v$ and 
$\|v\|'\ge\|u\|-\e$, and then a representation $u=\sum_{k=1}^mx^k\mmg z^k $ such that 
$\|u\|>\sum_{k=1}^m\|x^k\|\|z^k\|-\e$. Since evidently  $v=\sum_{k=1}^mx^k\odot z^k$, we have
 $\|v\|'\ge\|v\|_\odot-2\e$, and the reverse inequality follows. 

\smallskip
Now take $u\in X\mmg Z$. Let $(\dots,u^n,\dots)$ be the sequence $\odot_{X,Z}(u)$. We observe that 
$\rho_n^{X,Z}$ takes $u^n$ to $u_n:={\bf p}^n\cd u$; one can immediately check this on elementary 
tensors. It easily follows that $Ker(\s^{X,Z})=Ker(\odot_{X,Z})$. Since both $\s^{X,Z}$ and 
$\odot_{X,Z}$ are coisometries, there exists a unique isometric isomorphism $\iota_{X,Z}$, making 
the diagram (1.4) commutative. The rest is clear. $\tr$ 

\medskip
Thus, by virtue of Propositions 1.1--1.3, we have, for each $n$, a chain of isometric isomorphisms 
$$
X_n\mmp Z_n\stackrel{\rho_n^{X,Z}}{\longrightarrow}(X\mmg Z)_n\stackrel{\s_n^{X,Z}}{\longrightarrow}
(X\mmg Z)_n^{red}\stackrel{\iota_n^{X,Z}}{\longrightarrow}(X\odot Z)_n,\eqno(1.5)
$$
where the last map is the respective birestriction of $\iota^{X,Z}$. Denote by 
\\ $\vk_n^{X,Z}:X_n\mmp Z_n\to (X\odot Z)_n$ their composition. This is, of course, an isometric 
isomorphism of $c_0$-modules, well defined by taking $x\ot z$ to $x\odot z$. 

\medskip
{\bf Proposition 1.7}. {\it Suppose that at least one of modules $X$ and $Z$ is essential. Then the 
same is true for $X\mmg Z$, and, moreover, for every $u\in X\mmg Z$, we have } 
$$
u=\lim_{N\to\ii}P^N\cd u.\eqno(1.6)  
$$  

\smallskip
$\tl$ It follows from the equality (1.2), combined with the continuity of the operation ` $\mmg$ '. 
$\tr$ 

\medskip
From this, taking into account the diagram (1.4), we immediately obtain 

\medskip
{\bf Proposition 1.8}. {\it If at least one of modules $X$ and $Z$ is essential, then $\odot_{X,Z}$ 
is an isometric isomorphism of $c_0$-modules.}  $\tl\;\tr$

\medskip
The indicated assumption can not be omitted, even when both of modules are faithful: 

\medskip
{\bf Example 1.9}. Consider $X\!:=\!Z\!:=\!l_\ii$ with the coordinate-wise operations and uniform 
norm. 
Take the sequences $x:=(1,0,1,0,1,0,...)\in X$ and $z:=(0,1,0,1,0,1,...)\in Z$. Of course, we have 
$\odot_{X,Z}(x\mmg z)=0$. 

Now take two functionals $f,g:l_\ii\to\co$ of norm 1, such that $f(\xi)=g(\eta)=0$ for $\xi,\eta\in 
c_0$ and $f(x)=g(z)=1$; these are easily provided by the Hahn-Banach Theorem. Then the bilinear 
functional $f\times g:X\times Z\to\co:(\xi,\eta)\mt f(\xi)g(\eta)$ is obviously balanced and 
contractive. Therefore it gives rise to the contractive functional $f\mmg g:X\mmg Z\to\co$, well 
defined by $\xi\mmg\eta \mt f(\xi)g(\eta)$. Since $(f\mmg g)(x\mmg z)=1$, we have $x\mmg z\ne0$. 
Thus $\odot_{X,Z}$ is not injective. 

\bigskip 
Now suppose that we have three $c_0$-modules $X,Y$ and $Z$, so far arbitrary, and a bounded 
$c_0$-module morphism $\va:X\to Y$. The latter in an obvious way generates the sequence of its 
birestrictions $\va_n:X_n\to Y_n$. 

Consider the bounded morphism $\va\mmg\id: X\mmg Z\to Y\mmg Z$; we recall that it is well defined 
by $x\mmg z\mt\va(x)\mmg z$. Clearly, $\va\mmg\id$ maps $Ann(X\mmg Z)$ into $Ann(Y\mmg Z)$. It 
obviously follows that $\va\mmg\id$ gives rise to the bounded morphism $(\va\mmg\id)^{red}:   
(X\mmg Z)^{red}\to(Y\mmg Z)^{red}$, well defined by $(\va\mmg\id)^{red}(\s^{X,Z}(x\mmg z))= 
\s^{Y,Z}(\va(x)\mmg z); x\in X,z\in Z$. 

Combining this with Proposition 1.6, we obtain the commutative diagram 
$$
\xymatrix@C+20pt{X\mmg Z \ar[r]^{\odot_{X,Z}} \ar[d]_{\va\mmg\id}
&X\odot Z \ar[d]^{\va\odot\id\q}\\
Y\mmg Z \ar[r]^{\odot_{X,Z}} & {Y\odot Z} }\eqno(1.7)
$$
where $\va\odot\id$ is well defined by $x\odot z\mt\va(x)\odot z$. In other words, $\va\odot\id$ 
takes the sequence $(\dots,u_n,\dots); u_n\in X_n\mmp Z_n$ to the sequence 
$(\dots,(\va_n\ot\id)u_n,\dots)$. 

Note that we obviously have
$$
\|\va\odot\id\|\le\|\va\mmg\id\|\le\|\va\|\eqno(1.8)
$$

Being morphisms of $c_0$-modules, $\va\mmg\id$ and $\va\odot\id$ have well defined birestrictions 
$(\va\mmg\id)_n$ and $(\va\odot\id)_n$, respectively, for every $n$. Using the identifications, 
participating in the chain (1.5), for the pairs $(X,Z)$ and $(Y,Z)$, we easily obtain 

\medskip
{\bf Proposition 1.10}. {\it Both of $(\va\mmg\id)_n$ and $(\va\odot\id)_n$ are isometrically 
equivalent to the operator $\va_n\ot\id:X_n\mmp Z_n\to Y_n\mmp Z_n$. } $\tl\;\tr$ 

\medskip
Finally, Proposition 1.8 immediately implies 

\medskip
{\bf Proposition 1.11.} {\it Suppose that either both of $X$ and $Y$, or $Z$ are essential. Then 
the morphisms $\va\mmg\id$ and $\va\odot\id$ are isometrically equivalent}.  $\tl\;\tr$ 

\bigskip
{\centerline {\bf 2. Tensoring injective morphisms }} 

\bigskip
Let $X,Y,Z$ be normed $c_0$-modules, $\va:X\to Y$ a bounded morphism. Suppose that $\va$ is 
injective. When we can be sure that $\va\mmg\id$ is also injective ? (A kind of a `normed' version 
of an important typical question in pure algebra). 

If we ask the same about $\va\odot\id$, the situation is clear: 
 
\medskip
{\bf Proposition 2.1}. {\it Let $\va$ be injective. Then the same is true with $\va\odot\id$.} 

\smallskip
$\tl$ Together with $\va$, its birestrictions $\va_n$ are also injective. Then, for pure algebraic 
reasons, the same is true for operators $\va_n\ot\id: X_n\mmp Z_n\to Y_n\mmp Z_n$. It remains to 
recall the way $\va\odot\id$ acts. $\tr$ 

\medskip
From this we obtain 

\medskip
{\bf Proposition 2.2}. {\it Suppose that $X$ or $Z$ is essential. Then, if $\va$ is injective, 
then the same is true for $\va\mmg\id$}. 

\smallskip
$\tl$ By Propositions 1.8 and 2.1, both $\s^{X,Z}$ and $\iota^{X,Z}$ in the commutative diagram 
(1.7) are injective. The rest is clear. $\tr$  

\medskip
There is another kind of a condition, this time in terms of $\va$ itself, that gives the same 
result. Suppose that 
$\va$ is admissible, i.e. it has a left inverse bounded operator (not necessarily morphism of 
modules). Recall that the
Banach algebra $c_0$ is amenable, and hence every $c_0$-module, in particular, our $Z$, is flat. 
This means that for such a $\va$ the  morphism $\va\mmg\id$ is is not only injective, but 
topologically injective; see, e.g.,~\cite[Ch.VII]{he2}. (Actually, the cited book deals with the 
``completed'' theory, that is with Banach modules and completed module tensor products. But it is 
easy to observe that the indicated property of $\va\mmg\id$ is valid in the `non-completed' case as 
well). 

\medskip
However, if we have just an injective morphism between two normed $c_0$-modules, let them be even 
faithful, the situation is different: 

%


\medskip
{\bf Example 2.3}. Take $X:=Z:=l_\ii$ and set $Y:=c_0$. Consider a sequence 
$(\zeta_1,\zeta_2,...)\in c_0$ with non-zero terms and introduce $\va:X\to 
Y:(\xi_1,\xi_2,....)\mt(\zeta_1\xi_1,\zeta_2\xi_2,....)$. Of course, $\va$ is injective. At the 
same time, the lower horizontal arrow in (1.7) obviously depicts an injective map whereas the upper 
arrow, as we know from Example 1.9, does not. Therefore $\va\mmg\id$ can not be injective. 
 
\medskip
Of course, such a $\va$ is far from to be admissible. But what can happen in the ``intermediate'' 
case, when $\va$ is not bound to be admissible, but at least it is topologically injective? 

It is easy to show 
that $\va\mmg\id$ is not bound to be {\it topologically} injective. Moreover, as the related 
phenomenon,  in the `completed' theory such a morphism is not bound to be even injective (cf. the 
end of Introduction). 

But the present paper deals with the ``non-completed'' theory, and with a very specific base 
algebra. It turns out that in such a context we still have a positive result: 

\medskip 
{\bf  Theorem 2.4}. {\it Let $X,Y,Z$ be normed $c_0$-modules, and $\va:X\to Y$ be a topologically 
injective morphism. Then  $\va\mmg\id$ is injective.} 

\smallskip
$\tl$ Take $u\in X\mmg Z; u\ne0$; we want to show that $(\va\mmg\id)(u)$ is not 0. If 
$\odot_{X,Z}(u)\ne0$, that is $u\notin Ann(X\mmg Z)$, then the desired fact follows from 
Proposition 2.1, combined with the commutative diagram (1.7). Thus we have a right to assume that 
$u$ lives in $Ann(X\mmg Z)$. 

Consider the quotient maps $\tau_X:X\to X_{an}$ and $\tau_Z:Z\to Z_{an}$ (cf. Introduction) and 
set, for brevity, $\tau:=\tau_X\ot\tau_Z:X\ot Z\to X_{an}\ot Z_{an}$. Recall that $X\mmg Z$, by its 
definition, is a quotient space of $X\ot Z$ (actually, a quotient normed space of $X\mmp Z$) and 
denote by $\gamma$ the respective quotient map. It is easy to see that $Ker(\tau)$ is the 
  algebraic sum of $X_{es}\ot Z$ and $X\ot Z_{es}$. This obviously implies that
$$
\gamma(Ker(\tau))\su(X\mmg Z)_{es}.\eqno(2.1)
$$

Fix an arbitrary  $v\in X\ot Z$ with $\gamma(v)=u$ and set $w:=\tau(v)$.  We claim that $w\ne0$. 
 Indeed, in the opposite case we have, by (2.1), that $u\in(X\mmg Z)_{es}$ and hence, by (1.2), that
$u=\lim_{N\to\ii}P^N\cd u$. This, together with $u\in Ann(X\mmg Z)$, gives $u=0$, a contradiction. 

Thus $w$, being a non-zero vector in $X_{an}\ot Z_{an}$, can be represented as 
$w=\sum_{k=1}^n\tilde x_k\ot\tilde z_k; \tilde x_k\in X_{an}, \tilde z_k\in Z_{an}$, where $\tilde 
x_1\ne0$, and $\tilde z_k$ are linearly independent. 

\smallskip
Take an arbitrary $x_1\in X$ such that $\tau_X(x_1)=\tilde x_1$. Our next claim is that 
$\va(x_1)\notin Y_{es}$. Suppose the contrary. Then, by (1.2), we have 
$$
\va(x_1)=\lim_{N\to\ii}P^N\cd\va(x_1)=\lim_{N\to\ii}\va(P^N\cd x_1). 
$$
But this, since $\va$ is topologically injective, implies that $x_1=\lim_{N\to\ii}P^N\cd x_1$, that 
is $x_1\in X_{es}$. Hence we have $\tilde x_1=0$, a contradiction. 

This claim implies, by means of a standard corollary of the Hahn-Banach Theorem, that there exists 
a bounded functional $f:Y\to\co$ such that $f=0$ on $Y_{es}$, and $f(\va(x_1))=1$. The same 
corollary provides a bounded functional $\tilde g:Z_{an}\to\co$ such that $\tilde g(\tilde z_1)=1$ 
and $\tilde g(\tilde z_k)=0$ for $k=2,...,n$. Take an arbitrary $z_k\in Z$ with $\tau_Z(z_k)=\tilde 
z_k; k=1,...,n$ and consider the bounded functional $g:=\tilde g\tau_Z:Z\to\co$. Then we have, of 
course, that  $g(z_1)=1$ and $g(z_k)=0$ for $k=2,...,n$. 

Now introduce the bounded bilinear functional $f\times g:Y\times Z\to\co: \\ (y,z)\mt f(y)g(z)$. 
Since  $f=0$ on $Y_{es}$ and $g=0$ on $Z_{es}$, it is evidently balanced. Therefore it gives rise 
to the bounded linear functional, say $h:Y\mmg Z\to\co$, well defined by $h(y\mmg z)=f(y)g(z)$. 

We easily see that $h=0$ on $(Y\mmg Z)_{es}$. At the same time the element \\ $v-\sum_{k=1}^nx_k\ot 
z_k$ belongs to $Ker(\tau)$. Therefore we have, by (2.1), that \\ $u-\sum_{k=1}^nx_k\mmg 
z_k\in(X\mmg Z)_{es}$, and consequently $(\va\mmg\id)(u)-\sum_{k=1}^n\va(x_k)\mmg z_k$ lies in 
$(Y\mmg Z)_{es}$. Therefore $h(\va\mmg\id(u))=h(\sum_{k=1}^n\va(x_k)\mmg z_k)$, and the latter 
number is, of course, 1. It follows that $(\va\mmg\id)(u)\ne0$. $\tr$ 

\bigskip
{\centerline {\bf 3. Tensoring isometric morphisms }} 

\bigskip
In this section we shall deal with homogeneous $c_0$-modules, defined in Introduction. It is a 
rather large class of normed $c_0$-modules. In particular, we have 

\medskip
{\bf Proposition 3.1}. {\it Suppose that $X$ is an essential normed $c_0$-module, consisting of 
some complex-valued sequences and endowed with the coordinate-wise outer multiplication. Then $X$ 
is homogeneous.}
 
\smallskip
$\tl$ If $x,y\in X, x=(\dots,\lm_n,\dots), y=(\dots,\mu_n,\dots); \lm_n,\mu_n\in\co$, then 
 the equalities $\|x_n\|=\|y_n\|; n=1,2,...$ mean, of course, just that $|\lm_n|=|\mu_n|$.
 Therefore, for every $N\in{\B N}$ we have $P^N\cd x=\xi\cd P^N\cd y$ for some 
 $\xi=(\dots,\xi_n,\dots)\in c_0$ such that $|\xi_n|=1$ provided $n\le N$ and 
$\xi_n=0$ otherwise. It follows that $\|P^N\cd x\|\le\|P^N\cd y\|$, and similarly the reverse 
inequality is valid. But, since $X$ is essential, we can use (2). The rest is clear. $\tr$ 

\medskip
Note a useful

\medskip
{\bf Proposition 3.2.} {\it Let $X$ be a homogeneous $c_0$-module, $x\in X_{es}$ and $y\in X$. 
Suppose that $\|x_n\|\le\|y_n\|$ for all $n$. Then $\|x\|\le\|y\|$. } 

\smallskip
$\tl$ We have $\|x_n\|=\xi_n\|y_n\|$ for some $0\le\xi_n\le1; n=1,2,...$. Fix, for a moment, $N$, 
and consider $\xi:=(\xi_1,...,\xi_N,0,0,...)\in c_0$. Then, by homogeneity, we have $\|P^N\cd 
x\|=\|\xi P^N\cd y)\|\le\|y\|$. It remains to use (1.2) $\tr$ 

\medskip

\medskip
{\bf Proposition 3.3.} {\it Let $Z$ be a $c_0$-module. Assume that, for every essential homogeneous 
$c_0$-modules $X$ and $Y$ and an isometric morphism ${\textsf i}:X\to Y$ the morphism ${\textsf 
i}\mmg\id:X\mmg Z\to Y\mmg Z$ is also isometric. Then, for every $n=1,2,..$, the coordinate 
submodule $Z_n$ is, up to an isometric isomorphism of normed spaces, a dense subspace of 
$L_1(\Omega_n,\mu_n)$ for some measure space $(\Omega_n,\mu_n)$. } 

\smallskip
$\tl$ Suppose that, for a certain $n, Z_n$ does not satisfy the indicated condition. Then it easily 
follows from the criterion of Grothendieck~\cite[Thm. 1]{gro} that there are normed spaces $X,Y$ 
and an isometric operator ${\textsf i}:X\to Y$ such that the operator \\ ${\textsf i}\mmp\id:X\mmp 
Z_n\to Y\mmp Z_n$ fails to be an isometry. 

Set, for every $\xi=(\xi_1,...,\xi_n,...)\in c_0, x\in X, y\in Y$, $\xi\cd x:=\xi_nx$ and $\xi\cd 
y:=\xi_ny$. In this way we obviously make $X$ and $Y$ $c_0$-modules that are essential and 
homogeneous. Moreover, ${\textsf i}$ becomes a $c_0$-module morphism. Since $X$ and $Y$ are 
essential, it is sufficient, by virtue of Proposition 1.11, to show that the operator ${\textsf 
i}\odot\id:X\odot Z\to Y\odot Z$ is not an isometry. 

We see that, for $m\ne n$, we have $X_m=Y_m=0$. It easily follows that $X\odot Z=(X\odot Z)_n$ and 
$Y\odot Z=(Y\odot Z)_n$. Therefore the isometric isomorphisms $\vk_n^{X,Z}$ and $\vk_n^{Y,Z}$ (see 
Section 1) act between $X_n\mmp Z_n$ and $X\odot Z$, and, respectively, between $Y_n\mmp Z_n$ and 
$Y\odot Z$. Moreover, these isometric isomorphisms obviously implement an isometric equivalence of 
the operators ${\textsf i}\mmp\id$ and ${\textsf i}\odot\id$ (cf. (1.1)). Consequently, since the 
former of these two is not an isometry, the same is true for the latter. $\tr$ 

\medskip
Our principal aim is to show that the converse statement is valid. Actually, we shall prove a
slightly stronger assertion. 

\medskip
The main step in our proof is the following technical lemma. In what follows $S$ is an arbitrary 
homogeneous normed $c_0$-module with the following properties: 

(i) there exists a natural $N$ such that $S$, up to a linear isomorphism, is $\bigoplus_{n=1}^N 
S_n$. (In other words, for every $x\in S$ we have $P^N\cd x=x$). 

(ii) for every $n=1,...,N, S_n$ is a normed subspace of $L_1(\Omega_n,\mu_n)$ for some measure 
space $(\Omega_n,\mu_n)$, consisting of all step functions ( = linear combinations of 
characteristic functions of $\mu_n$-measurable subsets in $\Omega_n$). 

\medskip
{\bf Lemma 3.4}. {\it Let $X,Y$ be normed homogeneous $c_0$-modules and ${\textsf i}:X\to Y$ a 
morphism. Suppose we are given $u\in X\odot S$. Let $v:=({\textsf i}\odot\id_S)(u)\in Y\odot S$ be 
represented as $v=\sum_{k=1}^m y^k\odot g^k; y^k\in Y, g^k\in S$. Then for every $n=1,...,N$ there 
exist natural number $M, x^{kl}\in X_n$ and $g^{kl}\in S;  k=1,...,m, l=1,...,M$ such that for 
$$
y^{kl}:=y^k_1+y^k_2+\cdots+y^k_{n-1}+{\textsf i}_n(x^{kl}_n)+y^k_{n+1}+\cdots+y^N\eqno(3.1)
$$ 
we have 
$$
v=\sum_{k=1}^m\sum_{l=l}^M y^{kl}\odot g^{kl}\eqno(3.2)
$$
and }
$$
\sum_{k=1}^m\sum_{l=l}^M\|y^{kl}\|\|g^{kl}\|\le\sum_{k=1}^m\|y^k\|\|g^k\|.\eqno(3.3)
$$

\smallskip
$\tl$ Let $\sum_{s=1}^{m'}{'}x^s\odot f^s$ be an arbitrary representation of $u$. Remembering, what 
$S_n $ is, we can find $M\in{\B N}$ and a partition $\Om_n=\sqcup_{l=1}^M\Delta_l$, where 
$\Delta_l;l=1,...,N$ are $\mu_n$-measurable subsets of $\Om_n$ such that all $g^k_n,f^s_n$ are 
constant functions on each $\Delta_l$. In particular, for every $k=1,...,m, g^k_n$ has the form 
$\sum_{l=1}^M\lm^{kl}\chi_l$, where $\lm^{kl}\in\co$ and $\chi_l$ is the characteristic function of 
$\Delta_l$. 

Now for every $k=1,...,m,l=1,...,M$ we set
$$
 g^{kl}:=\frac{\|\lm^{kl}\chi_l\|}{\|g^k_n\|}g^k_1+\cdots+\frac{\|\lm^{kl}\chi_l\|}{\|g^k_n\|}g^k_{n-1}+
 \lm^{kl}\chi_l+\frac{\|\lm^{kl}\chi_l\|}{\|g^k_n\|}g^k_{n+1}+\cdots+\frac{\|\lm^{kl}\chi_l\|}
 {\|g^k_n\|} g^k_N\eqno(3.4)
$$
Since $\|\lm^{kl}\chi_l\|=\|\frac{\|\lm^{kl}\chi_l\|}{\|g^k_n\|}g^k_n\|$ and $S$ is homogeneous, we 
see that 
$$ 
\|g^{kl}\|=\|\frac{\|\lm^{kl}\chi_l\|}{\|g^k_n\|}g^k\| 
$$
for all $k,l$. But, living in $L_1(\cd)$, we have $\sum_{l=1}^M\|\lm^{kl}\chi_l\|=\|g^k_n\|$. 
Therefore for all $k$ we have 
 $\sum_{l=1}^M\frac{\|\lm^{kl}\chi_l\|}{\|g^k_n\|}=1$. Hence  
$g^k=\sum_{l=1}^M g^{kl}$ and 
$$
\|g^k\|=\sum_{l=1}^M\|\frac{\|\lm^{kl}\chi_l\|}{\|g^k_n\|}g^k\|=\sum_{l=1}^M\|g^{kl}\|.
$$ 
From this we have 
$$
v=\sum_{k=1}^m\sum_{l=l}^My^k\odot g^{kl}\q {\rm and} \q\sum_{k=1}^m\sum_{l=l}^M 
\|y^k\|\|g^{kl}\|=\sum_{k=1}^m\|y^k\|\|g^k\|.\eqno(3.5)
$$ 

\smallskip
Let us concentrate on $v_n$. It follows from (3.5) and (3.4) that 
$$
v_n=\sum_{k=1}^m\sum_{l=1}^My^k_n\otimes \lm^{kl}\chi_l=
\sum_{l=1}^M(\sum_{k=1}^m\lm^{kl}y^k_n)\ot\chi_l\eqno(3.6)
$$
But, as we remember, $v=({\textsf i}\odot\id_S)(u)$, and $u$ has the representation, indicated 
above. Therefore we have $v=\sum_{s=1}^{m'}{\textsf i}('x^s)\odot f^s$. Besides, by the choice  of 
$\Delta_l$, we have, for all $s$, that $f^s_n=\sum_{l=1}^M\nu^{sl}\chi_l$ for some 
$\nu^{sl}\in\co$.  Thus 
$$
v_n=\sum_{l=1}^M(\sum_{s=1}^m\nu^{sl}{\textsf i}_n('x^s))\ot\chi_l=\sum_{l=1}^M{\textsf i}_n(x^l)\ot\chi_l,
\eqno(3.7)
$$
where we set $x^l:=\sum_{s=1}^m\nu^{sl}('x^s)$.

 But $\chi_l; l=1,...,$  are linearly independent in $S_n$. Thus, comparing (3.7) and (3.6), we see that 
$$
\sum_{k=1}^m\lm^{kl}y^k_n={\textsf i}_n(x^l) \q {\rm  for\; all} \q l.\eqno(3.8)
$$ 

 Now introduce numbers 
 $$
 \al^{kl}:=(\lm^{kl})^{-1}\frac{\|\lm^{kl}y^k_n\|}{\sum_{t=1}^m\|\lm^{tl}y^t_n\|}\q {\rm provided} \q 
\lm^{kl}\ne0 \q {\rm and} \q \al^{kl}:=0 \q {\rm otherwise}.
$$ 
Finally, set $x^{kl}_n:=\al^{kl}_n x^{l}$.
 
\smallskip
Take $y^{kl}$ as in (3.1). Look at $v':=\sum_{k=1}^m\sum_{l=l}^M y^{kl}\odot g^{kl}$. By (3.1) and 
(3.5), $v'_{n'}=v_{n'}$ for all ${n'}\ne n$. As to $v'_n$, it is equal to 
$$
\sum_{k=1}^m\sum_{l=1}^M y^{kl}_n\ot g^{kl}_n=\sum_{l=1}^M \sum_{k=1}^m 
{\textsf i}_n(x^{kl})\ot \lm^{kl}\chi_l=
 \sum_{l=1}^M \sum_{k=1}^m {\textsf i}_n(\al^{kl}\lm^{kl}x^l_n)\ot\chi_l= 
 $$
 $$
 \sum_{l=1}^M \sum_{k=1}^m {\textsf i}_n\left(\frac{\|\lm^{kl} y^k_n\|}{\sum_{t=1}^m 
 \|\lm^{tl}y^t_n\|} x^l\right)\ot\chi_l= 
 \sum_{l=1}^M {\textsf i}_n(x^l)\ot\chi_l=\sum_{l=1}^M\lm^{kl} y^l_n\ot\chi_l,
 $$
 that is, by (3.8), to $v_n$. Thus $v$ and $v'$ have the same coordinates and hence, since $Y\odot S$ is 
 essential, they coincide. The equality (3.2) follows.

It remains to obtain (3.3). For this, we want to show that for all $l$ we have 
 $$
 \|{\textsf i}_n(x^{kl}_n)\|\le\|y^k_n\|.\eqno(3.9)
 $$
If $\al^{kl}=0$, this is immediate. Otherwise we have
$$
\|{\textsf i}_n(x^{kl}_n)\|=\|\al^{kl}{\textsf i}_n(x^{l})\|=
\left\|\frac{\|y^k_n\|}{\sum_{t=1}^m \|\lm^{tl}y^t_n\|}\left(\sum_{t=1}^m\lm^{tl}y^t_n\right)\right\|.
$$
and (3.9) follows from the triangle inequality for norms. 

Now it is time to use that $Y$ (not only $S$) is homogeneous. We have just shown that 
$\|y^{kl}_n\|\le\|y^k_n\|$, and, of course, we have $\|y^{kl}_{n'}\|=\|y^k_{n'}\|$ for all $n'\ne 
n$. Therefore, since all $y^{kl}$ belong to $Y_{es}$, Proposition 3.2 implies that 
 $$
 \|y^{kl}\|\le\|y^k\|.
 $$
Consequently, we have 
$$
\sum_{k=1}^m\sum_{l=l}^M\|y^{kl}\|\|g^{kl}\|\le\sum_{k=1}^m\sum_{l=l}^M\|y^{k}\|\|g^{kl}\|,
$$
and, because of (3.5), we are done. $\tr$ 

\medskip
{\bf Lemma 3.5}. {\it Let $X,Y,S$ be as in the previous lemma, and ${\textsf i}:X\to Y$ an 
isometric morphism. Then the morphism ${\textsf i}\mmg\id_S:X\mmg S\to Y\mmg S$ is also isometric.} 

\smallskip
$\tl$ Of course, $S$ is essential. Therefore, by virtue of Proposition 1.11, it is sufficient to 
prove that the morphism ${\textsf i}\odot\id_S:X\odot S\to Y\odot S$ is isometric. 

Fix an arbitrary $u\in X\odot S$ and set $v:=({\textsf i}\odot\id)(u)\in Y\odot S$. Our task is to 
show that $\|u\|=\|v\|$. 
 
 Take an arbitrary representation $v=\sum_{k=1}^m y^k\odot g^k; g^k\in S$. 
Set in the previous lemma $n:=1$. Getting rid of double sums, we can say that this lemma gives us a
representation 
$$
v=\sum_{k=1}^{m_1}y^{1k}\odot g^{1k},
$$
 where, for some $x^{1k}_1\in X_1, k=1,...,m_1$ and $y^{1k}_s, s=2,\dots,N$  we have
$$
y^{1k}={\textsf i}_1(x^{1k}_1)+y^{1k}_2+y^{1k}_3+\cdots+y^{1k}_N
$$
and
$$
\sum_{k=1}^{m_1}\|y^{1k}\|\|g^{1k}\|\le\sum_{k=1}^m\|y^k\|\|g^k\|.
$$

 Now apply Lemma 3.4 to the just obtained representation of $v$ and $n:=2$. 
 Looking at the form of the relevant $y^{kl}$  in the situation when the role of $y^k$ is  played by 
 $y^{1k}$ and again getting rid of double sums, we obtain a representation 
$$
v=\sum_{k=1}^{m_2} y^{2k}\odot g^{2k},
$$
 where, for some $x^{2k}_1\in X_1, x^{2k}\in X_2, k=1,...m_2$ and $y^{1k}_s, s=3,\dots,N$ , we have
$$
y^{2k}={\textsf i}_1(x^{2k}_1)+{\textsf i}_2(x^{2k}_2)+y^{3k}_3+\cdots+y^{2k}_N,
$$
and
$$
\sum_{k=1}^{m_1}\|y^{2k}\|\|g^{2k}\|\le\sum_{k=1}^{m_1}\|y^{1k}\|\|g^{1k}\| 
\q({\rm and\; hence} \q \le\sum_{k=1}^m\|y^k\|\|g^k\|).
$$

After this we apply Lemma 3.4 to this latter representation of $v$ and $n:=3$, and so on.  On the 
$N$th step, again (the last time) getting rid of double sums, we come to a representation of $v$ as 
$$
v=\sum_{k=1}^{m_N} y^{Nk}\odot g^{Nk},
$$
 where, for some $x^{Nk}_1\in X_1, x^{Nk}_2\in X_2,\dots,x^{Nk}_N\in X_N; k=1,...m_N$  we have
$$
y^{Nk}={\textsf i}_1(x^{2N}_1)+{\textsf i}_2(x^{Nk}_2)+\cdots+{\textsf i}_N(x^{Nk})
$$
and
$$
\sum_{k=1}^{m_N}\|y^{Nk}\|\|g^{2k}\|\le\sum_{k=1}^m\|y^k\|\|g^k\|.
$$

Now introduce $x^k:=x^{Nk}_1+\cdots+x^{Nk}_N\in X; k=1,\dots,m_N$. Obviously, $y^{Nk}={\textsf 
i}(x^k)$ and hence ${\textsf i}\odot\id_S(\sum_{k=1}^{m_N} x^{k}\odot g^{Nk})=v$. But ${\textsf 
i}\odot\id_S$ is injective (see Proposition 2.1). Therefore $\sum_{k=1}^{m_N} x^k\odot g^{Nk}$ is 
exactly $u$. Recalling that ${\textsf i}$ is isometric, we have 
$$
\|u\|\le\sum\|x^k\|\|g^{Nk}\|=\sum\|y^{Nk}\|\|g^{Nk}\|, 
$$
and hence
$$
\|u\|\le\sum_{k=1}^m\|y^k\|\|g^k\|.
$$
Taking the respective infimum in the expression (1.3) for the norm $\|\cd\|_\odot$, we have the 
estimate $\|u\|\le\|v\|$. Since, by (1.8), ${\textsf i}\odot\id$ is contractive, the desired 
equality follows. $\tr$ 

\medskip
{\bf Lemma 3.6}. {\it The assertion of the previous lemma remains true, if we replace the module 
$S$ by an arbitrary module $Z$ such that 

(i) there exists a natural $N$ such that $Z$ is linearly isomorphic to $\bigoplus_{n=1}^N Z_n$. 

(ii) 
for every $n=1,...,N, Z_n$ is, up to an isometric isomorphism, a dense normed subspace of 
$L_1(\Omega_n,\mu_n)$ for some measure space $(\Omega_n,\mu_n)$. }

\smallskip
$\tl$ Denote by $\bar Z$ and $\bar Z_n; n=1,\dots,N$ the completions of the $c_0$-modules $Z$ and 
$Z_n$, respectively. 

Take $z\in Z$. Obviously, we have 
$$
\max\{\|z_n\|; n=1,...,N\}\le\|z\|\le\sum_{n=1}^N\|z_n\|. 
$$
Therefore a sequence $z^m$ is a Cauchy  sequence in $Z$ if and only if for every $n=1,\dots,N$ the 
sequence $z^m_n$ is a Cauchy sequence in $Z_n$. It easily follows that $\bar Z$ is isometrically 
isomorphic to the algebraic direct sum $\bigoplus_{n=1}^N \bar Z_n$, endowed with the norm, well 
defined by $\|z\| =\lim_{m\to\ii}\|z^m\|$, where $z^m$ is an arbitrary sequence in $Z$ such that 
$\lim_{m\to\ii}z^m_n=z_n$ for every $n$. Obviously, $\bar Z_n$ is isometrically isomorphic to the 
space $L_1(\Om_n,\mu_n)$, mentioned in the formulation. It easily follows that $\bar Z$ contains 
the dense submodule $S$, satisfying the condition of Lemma 3.5. By virtue of that lemma, ${\textsf 
i}\mmg\id_S$ is an isometry. 

Therefore, by Proposition 1.2, the same is true for ${\textsf i}\mmg\id_{\bar Z}$, and this, in its 
turn, gives the desired property of ${\textsf i}\mmg\id_{Z}$. $\tr$ 

\medskip
{\bf Theorem 3.7.} {\it Let $Z$ be a  homogeneous $c_0$-module, satisfying the condition (ii) of 
the previous lemma. Further, let $X$ and $Y$ be two other homogeneous $c_0$-modules, ${\textsf 
i}:X\to Y$ an isometric morphism. Suppose that at least one of modules X and Z is essential. Then 
the morphism ${\textsf i}\mmg\id:X\mmg Z\to Y\mmg Z$ is also isometric.}   
  
$\tl$ Take $u\in X\mmg Z$. Our task is to show that $\|{\textsf i}\odot\id_Z(u)\|=\|u\|$. 

Fix, for a time, $N\in{\B N}$ and denote by $Z^{N}$ the submodule $\{P^N\cd u; u\in Z\}$ of $Z$. 
Consider the diagram 
$$
\xymatrix@C+20pt{X\mmg Z^{N} \ar[r]^{\id_X\mmg\, i} \ar[d]_{{\textsf 
i}'}
&X\mmg Z \ar[d]^{{\textsf i}'}\\
Y\mmg Z^{N} \ar[r]^{\id_Y\mmg\, i} & {Y\mmg Z} }
$$
where ${\textsf i}':={\textsf i}\mmg\id_Z$, and $i:Z^N\to Z$ is the natural embedding. By Lemma 
3.6, the left vertical arrow depicts an isometric morphism. Further, $\id_X\mmg i$ is contractive 
and has a contractive right inverse, namely $\id_X\mmg j$, where $j:Z\to Z^N$ acts as $z\mt P^N\cd 
z$. Therefore $\id_X\mmg i$ is an isometry, and the same is true with $\id_Y\mmg i$. 

For every $x\in X$ and $z\in Z$ we have $P^N\cd(x\mmg z)=x\mmg P^N\cd z$. From this, representing 
$u$ as a sum of elementary tensors, we see that $P^N\cd u=(\id_X\mmg i)(v)$ for some $v\in X\mmg 
Z^{(N)}$. Therefore, since our diagrum is obviously commutative and its three morphisms, mentioned 
above, are isometries, we have 
$$
\|({\textsf i}\mmg\id_Z)(P^N\cd u)\|=\|P^N\cd u\|. 
$$
Now observe that, by Proposition 1.7, we have $u=\lim_{N\to\ii}P^N\cd u$, and hence $\|({\textsf 
i}\odot\id_Z)(u)\|=\lim_{N\to\ii}\|({\textsf i}\odot\id_Z)(P^N\cd u)\|$. The rest is clear. $\tr$ 

\medskip
Combining this theorem with Proposition 3.3, we immediately obtain Theorem I, formulated in 
Introduction, with its mentioned corollaries for sequence modules and some other modules. 

From this theorem, in its turn, a Hahn-Banach type theorem, formulated in Introduction as Theorem 
II, easily follows. Indeed, it is a well known fact that, for a normed space $E$, its dual space is 
isometrically isomorphic to $L_\ii(\Om,\mu)$ for some measure space $(\Om,\mu)$ if and only if $E$ 
is isometrically isomorphic to a dense subspace of $L_1(\Om,\mu)$. (`If' part is the classics. To 
obtain the `only if' part we can recall, for example, that $L_\ii(\Om,\mu)$, being a von Neumann 
algebra, has only one, up to an isometric isomorphism, Banach predual space; cf., 
e.g.,~\cite[Cor.III.3.9]{tak}). Therefore, if we take this fact into account, Theorem II 
immediately follows from Theorem I, combined with Propositions 1.1 and 1.4.  

\newpage
{\centerline {\bf 4. A counter-example }} 
\bigskip

Here we want to show that the conditions in our main theorem, concerning the property of modules to 
be essential, can not be omitted, even within the class of faithful homogeneous modules. Namely, we 
shall show that the module $l_\ii$ (apparently the first faithful non-essential module that comes 
in mind), is not extremely flat with respect to the mentioned class. 

At first let us make some observations of general character. 

Let $X$ be a $c_0$-module. A subset $M$ of  ${\B N}$ is called a {\it support} of $X$, if we have 
$X_n=0$ for all $n\notin M$. 

\medskip
{\bf Lemma 4.1.}  {\it Let $X$ and $Z$ be two modules that have non-intersecting supports. Then for 
every $x\in X, x'\in X_{es}, z\in Z, z'\in Z_{es}$ we have $x'\mmg z=x\mmg z'=0$ in $X\mmg Z$.} 

\smallskip
$\tl$ By (1.2), we have 
$$x'\mmg z=\lim_{N\to\ii}P^N\cd x'\mmg z=\lim_{N\to\ii}\sum_{n=1}^N {\bf p}^n\cd x'\mmg {\bf p}^n\cd z.
$$ 
But the condition on supports implies that, for every $n$, either ${\bf p}^n\cd x'$ or ${\bf 
p}^n\cd z$ is 0. The rest is clear. $\tr$ 

\medskip
For $x\in X$, we shall denote by $\widetilde x$ the coset $x+X_{es}\in X_{an}$. 

\medskip
{\bf Proposition 4.2.} {\it Let $X$ and $Z$ be as before. Then there exists the isometric 
isomorphism of normed spaces $I_{X,Z}:X\mmg Z\to X_{an}\mmp Z_{an}$, well defined by \\ $x\mmg 
z\mt\widetilde x\ot\widetilde z$. } 

\smallskip
$\tl$ Consider the bilinear operator $X\times Z\to X_{an}\mmp Z_{an}:(x,z)\mt\widetilde 
x\ot\widetilde z$; it is obviously contractive and balanced. Therefore it gives rise to a 
contractive operator $I_{X,Z}$, well defined as it was indicated. 

Take $v\in X_{an}\mmp Z_{an}$, represented, say, as $\sum_{k=1}^n\widetilde x_k\ot\widetilde z_k; 
x_k\in X, z_k\in Z$. Then we have $v=I_{X,Z}(u)$, where $u=\sum_{k=1}^n x_k\mmg z_k$ with arbitrary 
$x_k,z_k$, taken in the respective cosets. Obviously, what we have to do is to show that 
$\|u\|\le\|v\|$. 

\smallskip
Take some $x_k'\in X_{es}, z_k'\in Z_{es}$. Lemma 4.1 implies that 

\smallskip
$u=\sum(x_k+x_k')\mmg(z_k+z_k')$. Therefore $\|u\|\le\sum_{k=1}^n\|x_k+x_k'\|\|z_k+z_k'\|$. Since 
$x_k',z_k'$ can be chosen in an arbitrary way, we have $\|u\|\le\sum_{k=1}^n\|\widetilde 
x_k\|\|\widetilde z_k\|$. Finally, since the taken representation of $v$ is also arbitrary, the 
very definition of the projective tensor norm gives the desired inequality. $\tr$ 

\medskip
 Now consider the normed quotient space (`ultraproduct') $l_\ii/c_0$. Since it is not isometrically 
 isomorphic to any space of the class $L_1(\Om,\mu)$, the theorem of Grothendieck, cited in Introduction, 
implies that there exist normed spaces $E$, $F$ and an isometric operator ${\widetilde i}:E\to F$ 
such that the operator 
$$
{\widetilde i}\mmp\id:E\mmp(l_\ii/c_0)\to F\mmp(l_\ii/c_0)
$$
 is not an isometry. Let us choose and fix these $E$, $F$ and ${\widetilde i}$. 

In what follows, we shall need, apart from the already used tensor product \\ ` $\mmp$ ', the  
non-completed {\it injective} tensor product of normed spaces and bounded operators, denoted by ` 
$\mmi$ ' (see, e.g.,~\cite[Ch.I.4]{def} or~\cite[Ch.3]{rya}). The injective tensor norm will be 
denoted by $\|\cd\|_i$. 

Consider the normed space $l_\ii\mmi E$. Evidently, it is a $c_0$-module with the outer 
multiplication well defined by $\xi\cd(\eta\ot x):=\xi\eta\ot x; \xi\in c_0, \eta\in l_\ii, x\in 
E$. 

This module is {\it contractive}: if $m_\xi:l_\ii\to l_\ii$ acts as $\eta\mt\xi\eta$, then, for 
every $u\in l_\ii\mmi E$, we have $\xi\cd u=(m_\xi\mmi\id_E)(u)$, and hence 
$$
\|\xi\cd u\|_i\le\|m_\xi\mmi\id_E\|\|u\|\le\|m_\xi\|\|\id_E\|\|u\|\le\|\xi\|\|u\|. 
$$

Besides, the introduced module is also {\it faithful}. Indeed, if $u\in l_\ii\mmi E$ is not 0, then 
it has a representation $u=\sum_{k=1}^n\eta^k\ot x^k$, where $x^k$ are linearly independent and 
$\eta^1\ne0$. Take $\xi\in c_0$ with $\xi\cd\eta^1\ne0$ and $g\in E^*$ with $g(x^1)\ne0$, 
$g(x^2)=...=g(x^n)=0$. Then $m_\xi\mmi g:l_\ii\mmi E\to l_\ii\mmi\co=l_\ii$ takes $\xi\cd u$ to 
$\xi\eta^1\mmi 1=\xi\eta^1$. Therefore $\xi\cd u\ne0$. 

Finally, the module $l_\ii\mmi E$ is {\it homogeneous}. This fact can be deduced from the known  
properties of the operation $\mmi C(\Om)$ (see, e.g., {\it idem}) and the identification of $l_\ii$ 
with $C(\beta{\B N})$. But we prefer to give a simpler proof. 

Obviously, it suffices to show that for $u\in l_\ii\mmi E$; $u=\sum_{k=1}^n\xi^k\ot x^k$ we have 
$$
\|u\|_i=\sup\{\|{\bf p}^n\cd u\|_i; n=1,2,...\}. 
$$

\smallskip
Take $f\in(l_\ii)^*$ and $g\in E^*$ with $\|f\|=\|g\|=1$. Then we have $(f\ot g)(u)=f(\eta^g)$, 
where $\eta^g:=\sum_{k=1}^ng(x^k)\xi^k$. Hence $|(f\ot 
g)(u)|\le\|\eta^g\|=\sup\{|(\eta^g)_n|;n=1,2,...\}$. But for every $n$ we have 
$$
|(\eta^g)_n|=\|\sum_{k=1}^n{\bf p}^n\xi^kg(x^k)\|=\|(\id\ot g)(\sum_{k=1}^n {\bf p}^n\xi^k\ot x^k\|=\|(\id\ot 
g)({\bf p}^n\cd u)\|\le\| {\bf p}^n\cd u)\|.
$$ 
Therefore the number $\|u\|_i$, which is, by definition, $\sup\{|(f\ot g)(u)|;f\in(l_\ii)^*, g\in 
E^*;\|f\|=\|g\|=1\}$, does not exceed $\sup\{\|{\bf p}^n\cd u\|_i; n=1,2,...\}$. Since the reverse 
inequality is obvious, we are done. 

\medskip
 In the same way we define the contractive faithful homogeneous $c_0$-module \\ $l_\ii\mmi F$. Finally, 
consider the operator ${\bf i}:=\id\mmi\widetilde i:l_\ii\mmi E\to l_\ii\mmi F$, which is evidently 
a morphism of $c_0$-modules. Because of the injective property of the operation ` $\mmi$ ' (see, 
e.g.,~\cite[Ch.I.4.3]{def} or~\cite[p. 47]{rya}), ${\bf i}$ is an isometry. 


From now on, it is convenient for us to use the notation $X$ for $l_\ii\mmi E$ and $Y$ for 
$l_\ii\mmi F$. 

\medskip
{\bf Theorem 4.3}. {\it The morphism ${\bf i}\mmg\id: X\mmg l_\ii\to Y\mmg l_\ii$ is not an 
isometry. As a corollary, the module $l_\ii$ is not extremely flat with respect to the class of all 
homogeneous normed $c_0$-modules}. 

\smallskip
$\tl$ We shall write $Z$ instead of $l_\ii$, and just $\id$ instead of $\id_Z$. Note that we have 
$Z_{an}=l_\ii/c_0$. 

Denote by $Z^{od}$ and $Z^{ev}$ the submodules of $Z$, consisting of  sequences  with the zero even 
terms and, respectively, zero odd terms. Besides, denote by $\id_{an}$ and $\id_\bullet$ the 
identity operators on $Z_{an}$ and, respectively, on $(Z^{ev})_{an}$. Our first claim is 

\medskip
1$^0$. The operator ${\widetilde i}\mmg\id_\bullet:E\mmp(Z^{ev})_{an}\to F\mmp(Z^{ev})_{an}$ is not 
an isometry. 

\smallskip Indeed, taking the sequence $(0,\xi_2,0,\xi_4,0,...)$ to $(\xi_2,\xi_4,...)$, we obtain 
isometric isomorphisms of normed spaces (by no means of modules) $j:Z^{ev}\to Z, \\ 
j_{es}:(Z^{ev})_{es}\to Z_{es}=c_0$ and, passing to respective cosets, $j_{an}:(Z^{ev})_{an}\to 
Z_{an}$. Then we easily see that the operators ${\widetilde i}\mmp\id_\bullet$ and ${\widetilde 
i}\mmp\id_{an}$ are isometrically equivalent. The rest is clear. 



\medskip
 From now on we shall use the brief notation $X^{od}$ for $Z^{od}\mmi E$, $Y^{od}$ for $Z^{od}\mmi F$, 
$\id^{od}$ for the identity operator on $Z^{od}$ and ${\bf i}^{od}$ for $\id^{od}\mmi\widetilde 
i:X^{od}\to Y^{od}$. Similarly to what was said about $X$ and $Y$, $X^{od}$ and $Y^{od}$ are 
contractive $c_0$-modules with respect to the same outer multiplication as for $X$ and $Y$ (cf. 
above), and ${\bf i}^{od}$ is an isometric morphism of $c_0$-modules. Besides, we introduce the 
operator 

\smallskip
${\bf i}_{an}:(X^{od})_{an}\to(Y^{od})_{an}$, which is well defined by taking a coset 
$x+(X^{od})_{es}$  to ${\bf i}^{od}(x)+(Y^{od})_{es}$. 

\medskip
Our next claim is 

\medskip
2$^0$. The operator  ${\bf i}_{an}\mmp\id_\bullet:(X^{od})_{an}\mmp(Z^{ev})_{an}\to 
(Y^{od})_{an}\mmp(Z^{ev})_{an}$ is not an isometry. 
 
\smallskip
Denote the sequence $(1,0,1,0,1,...)\in Z^{od}$ by $\widetilde1^{od}$. Consider the operator \\ 
${\bf s}_E:E\to(X^{od})_{an}$, taking a vector $x$ to the coset $(\widetilde1^{od}\mmi 
x)+(X^{od})_{es})$, and then  ${\bf s}_E\mmp\id_\bullet:E\mmp (Z^{ev})_{an}\to(X^{od})_{an}\mmp 
(Z^{ev})_{an}$. At first we shall show, as an intermediate step, that the latter operator is an 
isometry. 

For this aim, using the Hahn-Banach Theorem, introduce the functional \\ $h:Z^{od}\to\co$ of norm 
1, which takes the subspace $(Z^{od})_{es}=c_0\cap Z^{od}$ to 0 and takes $\widetilde1^{od}$ to 1. 
It gives rise to the operator ${\bf t}_E^0:=h\mmi\id_E:Z^{od}\mmi E\to\co\mmi E$, that is ${\bf 
t}_E^0:X^{od}\to E$. The latter evidently takes $(X^{od})_{es}$ to 0 and thererfore generates the  
operator ${\bf t}_E:=(X^{od})_{an}\to E$, well defined by taking the coset $u+(X^{od})_{es};u\in 
X^{od}$ to ${\bf t}_E^0(u)$. Since ${\bf s}_E$ and ${\bf t}_E$ are, of course, contractive, the 
same is true with ${\bf s}_E\mmp \id_\bullet$ and ${\bf t}_E\mmp\id_\bullet$. But the composition 
$({\bf t}_E\mmp\id_{an})({\bf s}_E\mmp\id_{an})$ is the identity operator on $E\mmp (Z^{ev})_{an}$. 
This implies that the former of these two is an isometry (and the latter is a coisometry). 

In a similar way, we introduce the operator \\
${\bf s}_F\mmp\id_\bullet:F\mmp (Z^{ev})_{an}\to(Y^{od})_{an}\mmp (Z^{ev})_{an}$ and show that it 
is also an isometry. 
Consequently, in the diagram 
$$
\xymatrix@C+20pt{E\mmp(Z^{ev})_{an} \ar[r]^{{\bf s}_E\mmp\id_\bullet} \ar[d]_{{\widetilde i}\mmp\id_\bullet}
&(X^{od})_{an}\mmp(Z^{ev})_{an} \ar[d]^{{\bf i}_{an}\mmp\id_\bullet}\\
F\mmp(Z^{ev})_{an} \ar[r]^{{\bf s}_F\mmp\id_\bullet} & {(Y^{od})_{an}\mmp(Z^{ev})_{an}} }
$$

\noindent the horizontal arrows depict isometries. Further, our diagram is obviously commutative. 
Thus it shows that the operator, depicted by the left vertical arrow, is isometrically equivalent 
to a birestriction of the operator, depicted by the right vertical arrow. But we already know that 
the former one is not an isometry. Therefore the same is true for the latter. 


\medskip
We turn to the next claim. 

\medskip
3$^0$. The morphism ${\bf i}^{od}\mmg\id^{ev}:X^{od}\mmg Z^{ev}\to Y^{od}\mmg Z^{ev}$ is not 
isometric. 

\smallskip
The set of odd natural numbers is the support of both $X^{od}$ and $Y^{od}$ whereas the set of even 
natural numbers is the support of $Z^{ev}$. Therefore Proposition 4.2 provides the isometric 
isomorphisms $I_{X^{od},Z^{ev}}:X^{od}\mmg Z^{ev}\to(X^{od})_{an}\mmp(Z^{ev})_{an}$ and 
$I_{Y^{od},Z^{ev}}:Y^{od}\mmg Z^{ev}\to({Y^{od})_{an}\mmp(Z^{ev})_{an}}$, well defined as it was 
indicated. Looking at the respective commutative diagram, we see that these isomorphisms implement 
the isometric equivalence between the operators ${\bf i}^{od}\mmg\id^{ev}$ and ${\bf 
i}_{an}\mmp\id_\bullet$. Thus the present claim follows from the previous one. 

 

\medskip
4$^0$. The end of the proof. 

\medskip 
Let $\rho^{od}:Z^{od}\to Z$ and $\rho^{ev}:Z^{ev}\to Z$ be the natural embeddings. Set 
$\rho^{od}_X:=\rho^{od}\mmi\id_E,\; \rho^{od}_Y:=\rho^{od}\mmi\id_F$; these maps are obviously 
morphisms of $c_0$-modules. Consider the diagram 
$$
\xymatrix@C+20pt{X^{od}\mmg Z^{ev} \ar[r]^{\rho^{od}_X\mmg\rho^{ev}} \ar[d]_{{\bf i}^{od}\mmg\id^{ev}}
&X\mmg Z \ar[d]^{{\bf i}\mmg\id}\\
Y^{od}\mmg Z^{ev} \ar[r]^{\rho^{od}_Y\mmg\rho^{ev}} & {Y\mmg Z} }
$$

Observe that its horizontal arrows depict isometries. Indeed, introduce the operators $\s^{od}:Z\to 
Z^{od}: (\xi_1,\xi_2,\xi_3,...)\mt(\xi_1,0,\xi_2,0,\xi_3,...)$, $\s^{ev}:Z\to Z^{ev}: 
(\xi_1,\xi_2,\xi_3,...)\mt(0,\xi_1,0,\xi_2,0,\xi_3,...)$ and set $\s^{od}_X:=\s^{od}\mmi\id_E:Z\mmi 
E\to Z^{od}\mmi E$. Obviously, the operator $\s^{od}_X\mmg\s^{ev}$ is  contractive, and the same is 
true with $\rho^{od}_X\mmg\rho^{ev}$. But the composition 
$$
(\s^{od}_X\mmg\s^{ev})(\rho^{od}_X\mmg\rho^{ev})=[(\s^{od}\rho^{od})\mmi\id_E]\mmg(\s^{ev}\rho^{ev})
$$
 is the identity operator on $X^{od}\mmg Z^{ev}$. This implies that the right factor, 
in our case $\rho^{od}_X\mmg\rho^{ev}$, is an isometry (whereas the left factor is a coisometry). 
Similarly, $\rho^{od}_Y\mmg\rho^{ev}$ is an isometry as well. 

Our diagram is clearly commutative, and, by the previous claim, its left vertical arrow does not 
depict an isometry. Hence the same is true with its right  vertical arrow (cf. the end of the proof 
of Claim 2). The rest is clear. $\tr$

\medskip
{\bf Remark}. The extreme flatness is a recent stronger version of a much older and more 
investigated notion of a strict flatness, that was mentioned in the introduction. We recall that 
the definition of a strictly flat module resembles that of an extremely flat module; one must only 
replace the  word `isometric' by `topologically injective' (see., e.g.,\cite{he3} ). 

The module $l_\ii$, as every normed module over the amenable algebra $c_0$, is (just) flat in the 
standard sense of~\cite{he2}~\cite{he3}\cite{run}. At the same time, by Theorem 4.3, it is not 
extremely flat. Here we want to note that one can show, using practically the same argument, as in 
the proof of the latter theorem, that it is not strictly flat as well. The only difference is that 
in the very beginning one must use a somewhat stronger property of $Z:=l_\ii/c_0$ than was employed 
before. Namely, there exist normed spaces $E,F$ and topologically injective operator ${\bf i}:E\to 
F$ such that ${\bf i}\mmp\id_Z$ is not topologicall injective. This is because $l_\ii/c_0$, being, 
in the terminology of~\cite{def}, an ${\cal L}^g_\ii$-space, can not be an ${\cal L}^g_1$-space 
(see Cor. 23.3(4) {\it idem}) This means, by Cor. 23.5(1) {\it idem}, that the operation $\mmp 
l_\ii/c_0$ `does not respect subspaces isomorphically' or, in our terminology, $l_\ii/c_0$ is not a 
strictly flat normed space ($\co$-module). The subsequent constructions and ``claims'' are, up to 
obvious modifications, the same.

\ed